\documentclass[11pt,a4paper]{article} 

\usepackage{setspace}
\usepackage{epsfig}
\usepackage{psfrag}
\usepackage[retainorgcmds]{IEEEtrantools} 
\usepackage{fullpage}
\usepackage{authblk}
\usepackage{latexsym}
\usepackage{amsmath,amssymb,amsfonts,amsthm}
\usepackage{graphicx}
\usepackage{natbib}
\usepackage[section]{placeins}
\newtheorem{definition}{Definition}
\newtheorem{theorem}{Theorem}

\newtheorem{corollary}{Corollary}

\makeatletter
\newcommand{\mbbR}{\mbox{$\mathbb{R}$}}
\newcommand{\mbbN}{\mbox{$\mathbb{N}$}}
\newcommand{\mbbP}{\mbox{$\mathbb{P}$}}
\newcommand{\RMSE}{\mbox{RMSE}}

\title{Dependence Properties of Multivariate Max-Stable Distributions}
\author[*]{Ioannis Papastathopoulos}
\author[$\dagger$]{Jonathan A. Tawn}
\affil[*]{\small School of Mathematics, University of Bristol}
\affil[$\dagger$]{\small Department of Mathematics and Statistics, Lancaster University}
\affil[$$]{\small i.papastathopoulos@bristol.ac.uk $\quad$ j.tawn@lancaster.ac.uk}

\date{}

\begin{document}
\maketitle
\begin{abstract}
  For an $m$-dimensional multivariate extreme value distribution there
  exist $2^{m}-1$ exponent measures which are linked and completely
  characterise the dependence of the distribution and all of its lower
  dimensional margins. In this paper we generalise the inequalities of
  \cite{schltawn02} for the sets of extremal coefficients and
  construct bounds that higher order exponent measures need to satisfy
  to be consistent with lower order exponent measures. Subsequently we
  construct nonparametric estimators of the exponent measures which
  impose, through a likelihood-based procedure, the new dependence
  constraints and provide an improvement on the unconstrained
  estimators.
\end{abstract}
\section{Introduction}
Max-stable distributions arise naturally from the study of limiting
distributions of appropriately scaled componentwise maxima of
independent and identically distributed random variables. Here and
throughout the vector algebra is to be interpreted as componentwise. A
vector random variable $Y=(Y_1,\hdots,Y_m)$ with unit Fr\'{e}chet
margins, i.e., $G_i(y):=\mbbP(Y_i < y)=\exp(-1/y)$, $y>0$, $i \in
M_m=\{1,\hdots,m\}$, is called max-stable if its distribution function
is max-stable, i.e., if
\begin{equation}
  G_{M_m}(y_{M_m}):=\mbbP\left(Y<y_{M_m}\right) = \exp\left\{- \int_{S_m} \max_{i\in
      M_m}\left(\frac{w_i}{y_i}\right)dH(w_1,\hdots,w_m)\right\},
  \label{eq:max_stable_distribution}
\end{equation}
where $y_{M_m} = (y_1,\hdots,y_m) \in \mbbR_+^m$, $S_m =
\left\{(w_{1},\hdots,w_m)\in \mbbR_{+}^{m}: \sum_{i=1}^m w_i =
  1\right\}$ is the $(m-1)$-dimensional unit simplex and $H$ is an
arbitrary finite measure that satisfies
\[
\int_{S_m}w_idH(w_1,\hdots,w_m)=1\quad \text{for any}\quad i\in M_m.
\]

\noindent The last condition is necessary for $G_{M_m}$ to have unit
Fr\'{e}chet margins and
representation~(\ref{eq:max_stable_distribution}) is due to
\cite{pick81}.\ There is no loss of generality in assuming unit
Fr\'{e}chet margins since our focus is placed on the dependence
structure of max-stable distributions, i.e., we are interested in the
copula function \citep{ne:99} which is invariant to strictly monotone
marginal transformations and in practice we can standardise random
variables to unit Fr\'{e}chet margins.  \newline

\noindent  The dependence properties of max-stable distributions have
received attention in the multivariate extreme value
literature. Dating back to \cite{sib60} and \cite{tiag62}, it has been
known that max-stable distributions are necessarily positively
quadrant dependent, i.e.,
\begin{equation}
G_{M_m}(y_{M_m}) \geq \prod_{i=1}^m G_i(y_i)\quad y_{M_m} \in \mbbR_+^m,
\label{eq:pqd}
\end{equation}
which implies that no pair of random variables can be negatively
dependent. Additionally, max-stable distributions satisfy even
stronger forms of dependence. \cite{marsolki83} show that
$\text{Cov}\{g(Y),h(Y)\}\geq 0$ for every pair of non-decreasing real
functions $g$ and $h$ on $\mbbR^m$, i.e., they are associated. For a
review of the dependence properties of max-stable distributions we
refer the reader to \cite{beiretal04} and the references
therein.\newline

\noindent Although all of the aforementioned properties exhibit
characteristics for the dependence structure of the class of
max-stable distributions, they are far too general to be either tested
or implemented in practice. In this paper, we introduce additional
constraints for the dependence structure that can be incorporated,
through a likelihood-based procedure, into the estimation of
max-stable distributions from observed componentwise maxima. The new
constraints are in essence the generalisation of the
\cite{schltawn02,schltawn03} inequalities for the extremal
coefficients which correspond to the dependence properties of
max-stable distributions for the special case of
$G_{M_m}(y,\hdots,y)$, $y>0$. As such, our notation and strategy are
influenced by the work of \cite{schltawn02,schltawn03}. The new
inequalities presented in this paper are related to the general case
of $G_{M_m}(y_{M_m})$, $y_{M_m}\in \mbbR_+^m$.  \newline

\noindent In Section~\ref{sec:dependence_properties} we introduce the
class of max-stable distributions along with the \cite{schltawn02}
inequalities for the extremal coefficients. Subsequently, we present
the general result of the paper that gives rise to inequalities for
the exponent measures. In Section~\ref{sec:inference} we consider the
\cite{halltajv00} nonparametric estimator for the exponent measure and
extend it, through a likelihood-based procedure, to satisfy the new
inequalities. Finally, in Section \ref{sec:simulation} a simulation
study is conducted to assess the performance of the constrained
estimator.
\newline

\section{Dependence Properties}
\label{sec:dependence_properties}
\subsection{Background}
\label{sec:back_and_not}
The class of max-stable distributions arises naturally from the study
of appropriately scaled component-wise maxima of random
variables. Consider a set of independent and identically distributed
random vectors $X^j=(X_1^{j},\hdots,X_m^{j})$, $j=1,\hdots,n$, with
unit Fr\'{e}chet margins. Under weak conditions \citep{resn87} it
follows that
\begin{equation}
  \lim_{n\rightarrow
    \infty}\mbbP\left(\bigcap_{i=1,\hdots,m}\left\{\max_{j=1,\hdots,n}X_{i}^{j}/n
      < y_i\right\}\right) = G_{M_m}(y_{M_m}), \quad y_{M_m}\in
  \mbbR_+^m.
\label{eq:limiting_law}
\end{equation}
The distribution function $G_{M_m}$ can be completely characterised by
the following representations
\begin{IEEEeqnarray}{rCl}
  V_{M_m}(y_{M_m})=-\log G_{M_m}(y)&=&\int_{S_m}\max_{i\in
    M_m}\left(\frac{w_i}{y_i}\right)dH(w_1,\hdots,w_m),
  \label{eq:Gspmeas}\\\nonumber\\
  & = &\left\{\sum_{i=1}^{m} 1/y_i\right\}
  A_{M_m}\left(\frac{1/y_1}{\sum_{i=1}^{m}
      1/y_i},\hdots,\frac{1/y_m}{\sum_{i=1}^{m} 1/y_i}\right),
  \label{eq:pickands_A}
\end{IEEEeqnarray}
where the function $V_{M_m}$ is known as the exponent measure of the
multivariate extreme value distribution $G_{M_m}$ and $A_{M_m}$,
called the \emph{Pickands' dependence function}, is a convex function
that satisfies $\max\{w_1,\hdots,w_m\}\leq A_{M_m}(w_1,\hdots,w_m)
\leq 1$, $(w_1,\hdots,w_m) \in S_m$. This condition implies that
$A_{M_m}(e_j)=1$, $j \in M_m$, where $e_j$ is the $j$-th unit vector
in $\mbbR^m$. \newline

\noindent Let $C_m=2^{M_m}\setminus \{\emptyset\}$ and denote also by
$y_B=\{y_i: i\in B\}$ for $B \in C_m$. Then we can define $2^m - 1$
exponent measures for an $m$-dimensional max-stable random vector $Y$,
where each one characterises completely the distribution function of a
marginal random variable $Y_B$ of $Y$, i.e.,
\[
V_B(y_B) = -\log\left\{\mbbP(Y_B < y_B)\right\} =
-\log\left\{\lim_{y_{M_m\setminus B}\rightarrow \boldsymbol{\infty}
  }G_{M_m}(y_{M_m})\right\}, \quad B \in C_m.
\]
The set of exponent measures $\{V_B: B \in C_m\}$ describes completely
the dependence structure of a max-stable distribution given by
equation~(\ref{eq:max_stable_distribution}) and all of its lower
dimensional margins. It is also trivial to see that with each exponent
measure $V_B$ there is an associated Pickands' dependence function
$A_B$. Additionally, $V_B$, $B \in C_m$, is homogeneous of order $-1$,
i.e.,\ $V_B\left(y,\hdots,y\right)=y^{-1} V_B\left(1,\hdots,1\right)$,
$y>0$.\newline 

\noindent The importance of the homogeneity property is mostly
illustrated through one widely used measure of extremal dependence for
the variables indexed by a set $B\in C_m$.\ More specifically, the
quantity defined by
\begin{equation}
\theta_B = V_B(1,\hdots,1) = \int_{S_m}\max_{i \in B}w_i
dH(w_1,\hdots,w_m), \quad 1\leq \theta_B \leq |B|,
\label{eq:ext_coeff}
\end{equation}
describes the effective number of independent variables in the set $B$
and arises naturally from the distribution of the maximum of all the
variables indexed by the set $B$, i.e.,
\begin{equation}
\mbbP\Big\{\max_{i\in B}Y_i < y\Big\}= \mbbP\left\{Y_i <
  y\right\}^{\theta_B}, \quad y>0.
\label{eq:ext_coeff_power}
\end{equation}
The measure $\theta_B$ is termed the extremal coefficient and complete
dependence and independence corresponds to $\theta_B=1$ and
$\theta_B=|B|$ respectively.\ Also, from
expression~(\ref{eq:ext_coeff_power}) it follows trivially that
$\theta_{B}=1$ for any $B \in C_m$ with $|B|=1$. Due to its simple
interpretation, the set of extremal coefficients $\{\theta_B: B\in
C\}$ has been used as a dependence measure in various applications
\citep{tawn90,schltawn03}.

\subsection{\cite{schltawn02,schltawn03} inequalities for the extremal
  coefficients}
\label{sec:schltawn02}
\cite{schltawn02,schltawn03} constructed bounds for the set of
extremal coefficients $\{\theta_{B}:B \in C_m\}$ of max-stable
distributions that characterise the dependence structure for the
special case of $G_{M_m}(y,\hdots,y)$, $y>0$. Here we use the
terminology of \cite{schltawn02} and for non-empty distinct subsets
$B_1,\hdots,B_s$ of $M_m$, $s \in \mbbN$, we refer to the set of
extremal coefficients $\{\theta_{B_1},\hdots,\theta_{B_s}\}$ as
complete and consistent if $s=2^m - 1$ and $\theta_{B_i}$ is given by
expression~(\ref{eq:ext_coeff}), respectively. Their main result is
given in the following theorem.

\begin{theorem}[\cite{schltawn02} Corollary 5]
  \label{th:schltawn02}
  A complete set of extremal coefficients $\{\theta_B: B \in
  C_m\}$, where $M_m$ is a finite set of indices, is consistent if
  and only if 
  \begin{equation}
    \sum_{B \in C_m, B \supseteq M_m\setminus L} (-1)^{|B \cap L|+1} \theta_B \geq 0, \quad \text{for all $L\in C_m$}.
    \label{eq:schltawn02_inequalities}
  \end{equation}
\end{theorem}

\noindent Theorem~\ref{th:schltawn02} yields bounds that higher order
extremal coefficients need to satisfy to be consistent with lower
order extremal coefficients. For example consider the
inequalities~(\ref{eq:schltawn02_inequalities}) for the cases $m=2$
and $m=3$ and let for ease of notation $\theta_{\{i,j\}}$ and
$\theta_{\{i,j,k\}}$ be $\theta_{ij}$ and $\theta_{ijk}$, for $i,j,k
\in M_m$ and $i \neq j \neq k$. These are respectively
\begin{IEEEeqnarray*}{CC}
  &1 \leq \theta_{12},\theta_{13},\theta_{23} \leq 2\\ & \text{and}\\
  &\max\left\{\theta_{12},\theta_{13},\theta_{23},\theta_{12}+\theta_{13}+\theta_{23}
    - 3\right\} \leq \theta_{123} \leq
  \min\left\{\theta_{12}+\theta_{13}-1,\theta_{12}+\theta_{23}-1,\theta_{13}+\theta_{23}-1\right\}.
\end{IEEEeqnarray*}
The first set of inequalities represents the well known bounds of the
extremal coefficients that come from the positively quadrant
dependence property~(\ref{eq:pqd}) of max-stable
distributions. However, the second set of inequalities gives tighter
bounds for the higher order extremal coefficient $\theta_{123}$. This
can be seen easily since the combined inequalities for the cases $m=2$
and $m=3$ reduce to $1\leq \theta_{123}\leq 3$.

\subsection{Inequalities for the exponent measures of max-stable
  distributions}
\label{sec:inequalities}
It transpires that similar inequalities as with those in
expression~(\ref{eq:schltawn02_inequalities}) can be obtained for the
exponent measures $\{V_B:B\in C_m\}$ of max-stable
distributions. Analogously with the terminology for the extremal
coefficients in Section~\ref{sec:schltawn02} we introduce the
following definition.
\begin{definition}
  \label{def:consistent_V}
  Let $s$ be an integer, for $i=1,\hdots,s$ $B_i$ are distinct
  non-empty subsets of $M_m=\{1,\hdots,m\}$ and
  $y_{M_m}=(y_1,\hdots,y_m) \in \mbbR_+^{m}$. An ensemble
  $\left\{V_{B_1}(y_{B_1}),\hdots,V_{B_s}(y_{B_s})\right\}$ of
  exponent measures, where $y_{B_i}=\{y_j:j\in B_i\}$, is called
  consistent if
  \[
  V_{B_i}(y_{B_i}) =\int_{S_m} \max_{j\in B_i}
  \left(\frac{w_j}{y_j}\right) dH(w_1,\hdots,w_m),
  \]
  for $i=1,\hdots,s$ and $H$ is an arbitrary finite measure that
  satisfies $\int_{S_m}w_idH(w_1,\hdots,w_m)=1$ for any $i\in M_m$.\
\end{definition}
\noindent If $s=2^m-1$ then the set of exponent measures is called
complete. The following theorem provides a new representation of the
exponent measures of multivariate extreme-value distributions in terms
of non-negative and uniquely defined real functions.
\begin{theorem}
  \label{th:theorem1}
  Let $\{V_B:B\in C_m\}$ be a complete and consistent set of exponent
  measures. Then, there exist $2^m-1$ non-negative functions
  $d_L:\mbbR_+^{m}\rightarrow \mbbR_+$, $L\in C_m$, such that, for any
  $B \in C_m$
  \begin{equation}
    V_{B}(y_B) = \sum_{L\in M_m, L\cap B \neq \emptyset} d_L\left(y_{M_m}\right),
    \label{eq:drepresentation}
  \end{equation}
  and the functions $d_L$ are uniquely given by
  \begin{equation}
    d_L\left(y_{M_m}\right) = \sum_{B\in C_m, B \supseteq M_m \setminus L} 
    \left(-1\right)^{|B\cap L|+1} \int_{S_m} \max_{j\in B} 
    \left(\frac{w_j}{y_j}\right) dH(w_1,\hdots,w_m).
    \label{eq:d_L}
  \end{equation}
\end{theorem}
\subsubsection*{Proof}
  The proof of equation~(\ref{eq:drepresentation}) of
  Theorem~\ref{th:theorem1} follows along the lines of
  \cite{schltawn02} proof of Theorem~5 for the simpler case of the
  extremal coefficients by replacing the constants $\alpha_{k}^{i}(n)$
  of \cite{deheuv83} representation of max-stable distributions with
  $\alpha_{k}^{i}(n)/y_i$, $i\in M_m$, $k\in
  \mathbb{Z}$. Equation~(\ref{eq:d_L}) is the M\"{o}bius inversion of
  equation~(\ref{eq:drepresentation}).\qedsymbol\newline
\noindent The characterisation of a consistent set of exponent
measures is obtained from the following corollary.
\begin{corollary}
  \label{cor:nonneg}
  A complete set of exponent measures $\{V_B: B\in C_m\}$ is
  consistent if and only if 
  \begin{equation}
    \label{eq:ineqgenerator}
    \sum_{B\in C_m, B \supseteq M_m \setminus L} \left(-1\right)^{|B\cap
      L|+1} \int_{S_m} \max_{j\in B} \left(\frac{w_j}{y_j}\right)
    dH(w_1,\hdots,w_m) \geq 0,
  \end{equation}
  for all $y_{M_m}\in \mbbR_+^{m}$ and $L\in C_m$.
\end{corollary}

\section{Inference}
\label{sec:inference}
\subsection{The \cite{halltajv00} estimator of the exponent measure}
\label{sec:ex_ests}
The fundamental premise in all statistical extreme value modelling is
that the observed extremes of a stochastic process are well modelled
by the limiting theoretical extreme-value distributions. Let for
example $X^j=(X_1^{j},\hdots,X_m^{j})$, $j=1,\hdots,N$, be a set of
independent and identically distributed $m$-dimensional random vectors
with unit Fr\'{e}chet margins. Here and throughout we assume that the
normalised componentwise block maxima
\[
Y^j:= \bigvee_{r = (j-1) d + 1}^{j d} \frac{X^r}{d}, \quad j =
1,\hdots,n,
\]
where $n d = N$, follow exactly the law $G_{M_m}$ of the limiting
expression~(\ref{eq:limiting_law}).\newline

\noindent Let now $w_{B}\in S_{|B|} = \left\{w_{B}\in \mbbR^{|B|}_+ :
  \sum_{i\in B}w_{B,i}=1\right\}$, $B \in C_m$, and define $Z_{B}^j
= w_{B} Y_{B}^j$, for $j=1,\hdots,n$. It then follows that the
cumulative distribution function of $\max_{i \in B} Z_{i}^j$ is
Fr\'{e}chet with scale parameter equal to the Pickands' dependence
function $A_{B}(w_{B})$ of $G_{B}$, i.e.,
\[
\mbbP\left\{\max_{i \in B} Z_{i}^{j} < y\right\} = \exp\left\{-
  \frac{A_{B}(w_{B})}{y}\right\}, \quad y>0.
\]
A natural consistent estimator of $A_{B}$ then is the
\cite{halltajv00} corrected version of Pickands' estimator
\citep{pick81} which maximises the likelihood
\begin{equation}
  \ell_B\left\{A_{B}(w_B)\right\} = n
  \log\left\{A_{B}(w_B)\right\} - 2 \sum_{j=1}^{n} \log W_{B}^j -
  A_{B}(w_B)\sum_{j=1}^{n} \frac{1}{W_{B}^j},
\label{eq:llik}
\end{equation}
where $W_{B}^j = \max_{i\in B} \left\{w_{B,i} Y_{i}^j
  \left[\sum_{j=1}^n (1/Y_{i}^j)/n\right]\right\}$, $j=1,\hdots,n$,
is the \cite{halltajv00} correction which ensures that $\max
w_B\leq\hat{A}_{B}(w_B)$, for all $w_B\in S_{|B|}$, as well as
$\hat{A}_B(e_j) = 1$, for any $j \in M_m$, where $e_j$ is the $j$-th
unit vector in $\mbbR^m$. The maximum likelihood estimator is given by
$\hat{A}_{B}(w_B) = \left\{ n^{-1 }\sum_{j=1}^{n}
  (1/W_B^j)\right\}^{-1}$ which is subsequently corrected by
\[
\tilde{A}_{B}(w_B)=\min\left\{\hat{A}_{B}(w_B),1\right\}
\]
to satisfy $\tilde{A}_{B}(w_B)\leq 1$, for all $w_{B} \in S_{|B|}$. On
combining the estimator $\tilde{A}_{B}$ with
equation~(\ref{eq:pickands_A}), the following consistent estimator of
the exponent measure $V_{B}$ is obtained,
\begin{equation}
  \tilde{V}_{B}(y_{B}) = \left\{\sum_{i \in B} 1/y_i\right\}
  \tilde{A}_{B}\left(\frac{1/y_B}{\sum_{i\in B} 1/y_{i}}\right), \quad y_B\in \mbbR_+^{|B|},\quad B \in C_m.
  \label{eq:V_est}
\end{equation}
Other types of estimators exist in the literature such as the
non-parametric estimators proposed by \cite{deuh91} and
\cite{capetal97} for the bivariate case. \cite{zhangetal08} gives a
detailed overview of the existing estimators and extends them to the
multivariate case. In this paper though we use the \cite{halltajv00}
estimator since it arises as the maximum of a log-likelihood function
based on which the new inequalities of Section~\ref{sec:inequalities}
can be imposed.



\subsection{Constrained estimators}
\label{sec:constr_estimators}   
It transpires that the aforementioned nonparametric estimators of the
exponent measures do not necessarily ensure that the resulting
estimated set of exponent measures satisfy
inequalities~(\ref{eq:ineqgenerator}). The focus here is placed on
incorporating these additional constraints in the estimation procedure
so that the resulting complete set of estimated exponent measures
$V_B(y_B)$, $B\in C_m$, is consistent in the sense of
Definition~\ref{def:consistent_V} for fixed $y_{M_m}\in \mbbR_+^m$. To
incorporate the inequalities we construct similarly with
\cite{schltawn03} a joint log-likelihood function $\ell$ of
$\left\{A_B; B\in C_m\right\}$ by falsely assuming independence
between the observations for all different $B$ to give the
pseudo-log-likelihood
\begin{equation}
\ell\left(\left\{A_B\left(\frac{1/y_B}{\sum_{i\in B} 1/y_{i}}\right):
    B\in C_m\right\}\right) = \sum_{B \in C_m, |B|\geq 2 } \ell_{B}
\left\{ A_{B}\left(\frac{1/y_B}{\sum_{i\in B} 1/y_{i}}\right)
\right\}.
\label{eq:pseudo_likelihood}
\end{equation}
The maximum pseudo-likelihood estimators are consistent
\citep{lianself96} and the constrained estimators are obtained by
maximising the pseudo-log-likelihood~(\ref{eq:pseudo_likelihood})
subject to
\[
\sum_{B\in C_m, B \supseteq M_m \setminus L} \left(-1\right)^{|B\cap
  L|+1} \left\{\sum_{i \in B}1/y_i\right\}
A_{B}\left(\frac{1/y_B}{\sum_{i \in B} 1/y_i}\right) \geq 0,\quad \text{for
  all $L\in C_m$}\quad
\]

\noindent and

\[  
A_{B}\left(\frac{1/y_B}{\sum_{i \in B} 1/y_i}\right)\leq 1, \quad
\text{for all $B\in C_m$}.
\]
The resulting constrained estimators are denoted by $\tilde{A}_B^c$
which in turn yield the estimators $\tilde{V}_B^c$ as in
equation~(\ref{eq:V_est}). The joint estimation of the exponent
measures ensures that all estimators are self-consistent. Note that
the resulting estimates of lower order exponent measures are affected
by higher order measures, i.e., estimates of $V_{B_0}(y_{B_0})$ are
affected by estimates of $V_{B_1}(y_{B_1})$, where $B_0 \subset
B_1$. The major benefit of this feature is that this guarantees the
existence of higher order measures which are self-consistent with the
lower order measures.\newline

\noindent An alternative way of obtaining a set of estimated exponent
measures is via sequential estimation, i.e., the lower order exponent
measures are estimated firstly and then are used as constraints in the
estimation of the higher order exponent measures, see also
\cite{schltawn03}. Although this method is faster than the joint
optimization problem described by
equation~(\ref{eq:pseudo_likelihood}), it does not have the desirable
feature described above.

\section{Simulation Study}
\label{sec:simulation}

\subsection{Design}
\label{sec:design}

We illustrate the impact of constraining the \cite{halltajv00}
estimators to satisfy the new inequalities~(\ref{eq:ineqgenerator})
over the unconstrained estimators of the set of exponent measures
using simulated data from a 3-dimensional max-stable distribution,
i.e., the extreme value logistic distribution with dependence
parameter $\alpha \in (0,1]$ and set of exponent measures given by
\begin{equation}
  \left\{  V_{B}\left(y_B\right)= \left(\sum_{i\in B} y_i^{-1/\alpha}\right)^\alpha: \quad B\in C_3\right\},\quad y_B\in \mbbR_+^{|B|}.
  \label{eq:logistic}
\end{equation}
The values $\alpha=1$ and $\alpha=0$, taken as $\alpha\rightarrow 0$,
correspond to independence and complete dependence,
respectively. \newline

\noindent All comparisons are based on the root mean square error
(\RMSE) performance of the exponent measure estimators for a range of
dependence parameters $\alpha$ and a cube grid of values, say
$\mathbb{L}^3\subseteq \mbbR_+^3$, for $y_{M_3}$. Specifically, the
values chosen for the dependence parameter and the sample size are
$\alpha \in \{0.2, 0.5, 0.8\}$ and $n=50$, respectively. Results from
larger sample sizes are not reported in the paper since they are
unrealistic for applications and also, the efficiency of the
estimators $\tilde{V}_{B}$ and $\tilde{V}_{B}^c$ is similar, a fact
that comes from the consistency property of the \cite{halltajv00}
estimator. The set $\mathbb{L}$ was chosen to be the discrete set
$\{x_{p_1},\hdots,x_{p_7}\}$ with $x_p$ denoting the $p$-th quantile
of the unit Fr\'{e}chet distribution. We chose $p_1 = 0.05$,
$p_7=0.95$ and step size $p_j - p_{j-1} = 0.15$. The Monte Carlo size
used to compute estimates of the RMSE is 500.\newline

\noindent To obtain an aggregated measure of performance, we also
report the Monte Carlo estimates of the integrated square deviation of
the estimators from the theoretical function, i.e.,
\begin{equation}
\tilde{T}_{B} = \int_{C\left(\mathbb{L}^{|B|}\right)}\left\{\tilde{V}_B(y_{B}) -
  V_B(y_{B})\right\}^2d y_{B},\quad \text{for all}\quad B \in C_3,
\label{eq:T}
\end{equation}
\noindent where $C(\mathbb{L}^{|B|})$ is the smallest $|B|$-hypercube
that contains the set $\mathbb{L}^{|B|}$. The integral in
expression~(\ref{eq:T}) is approximated in each Monte Carlo iteration
by the quadrature mid-point numerical integration technique on the
grid $\mathbb{L}^{|B|}$ and the measure $\tilde{T}_B^c$ is defined
analogously by replacing $\tilde{V}_B$ in expression~(\ref{eq:T}) with
$\tilde{V}_B^c$.
\subsection{Results}
\label{sec:results}
\begin{figure}[htbp!]
  \centering
  \includegraphics[scale=.92,trim= 5 5 5 5]{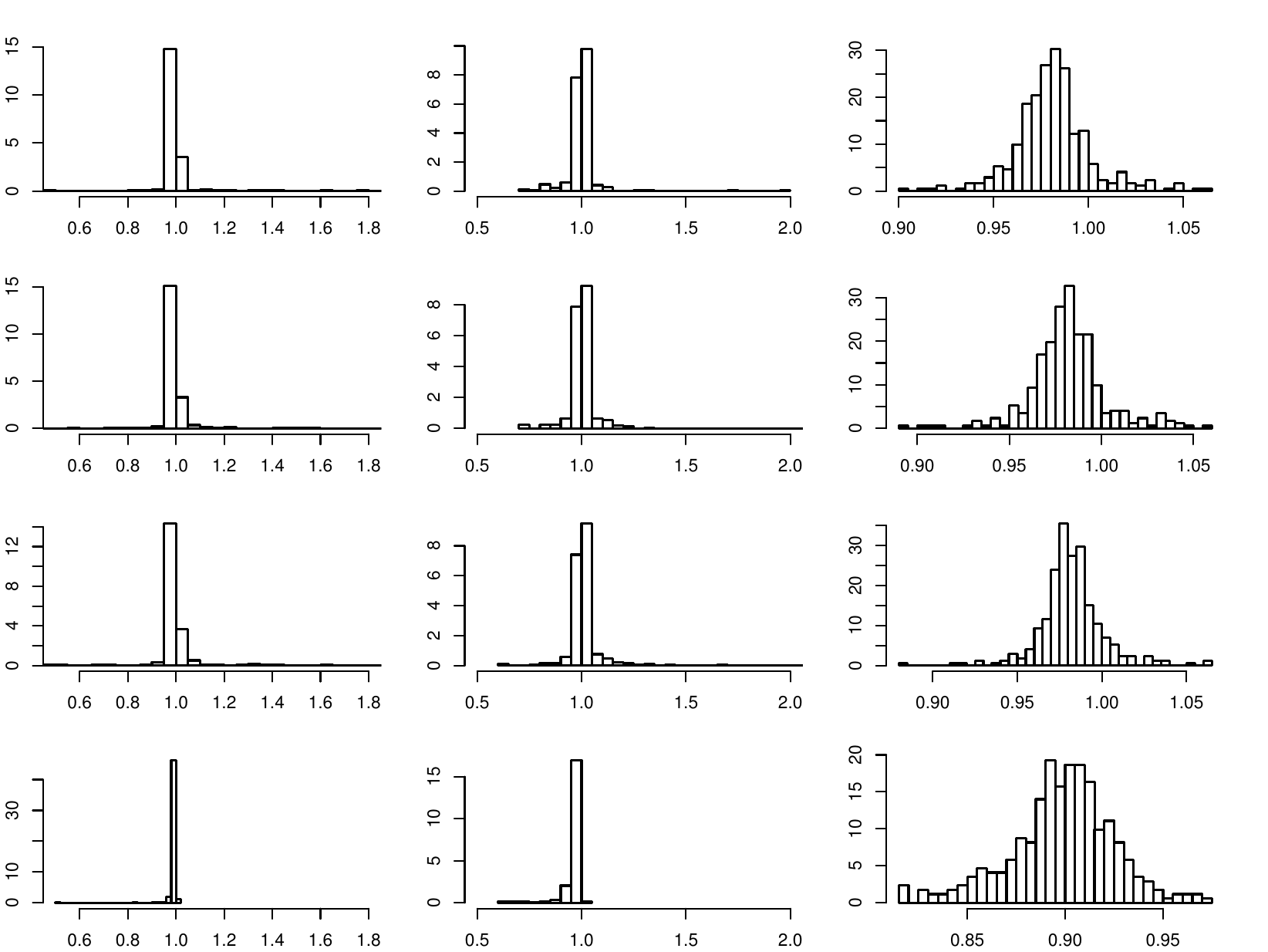}
  \caption{Histograms of the ratio of Monte Carlo estimates of RMSEs
    between the constrained and unconstrained estimators of the
    exponent measures $V_{12},V_{13},V_{23}$ and $V_{123}$ (top to
    bottom) for the extreme value logistic distribution with $\alpha
    = 0.2$ (left) $\alpha=0.5$ (centre) and $\alpha=0.8$
    (right).}
\label{fig:ratio_rmse}
\end{figure}

\noindent Figure~\ref{fig:ratio_rmse} shows the histograms of the
ratio of RMSEs between $\tilde{V}_B^c$ and $\tilde{V}_B$, $B \in C_3$,
for all grid points in $\mbbR_+^3$ for the extreme value logistic
distribution with $\alpha = 0.2,0.5$ and $0.8$. The figures indicate
either similar or better performance of the constrained estimators.

\noindent In particular, for the $\alpha=0.8$ case, the constrained
estimators are more efficient than the unconstrained estimators
especially for the higher order exponent measure $V_{123}$ and
improvement in RMSE, although lower in magnitude, can be also seen in
the bivariate exponent measures $V_{B}$, $B\in C_3 \setminus
M_3$. Also, the percentage of Monte Carlo samples where the
constrained estimates changed with respect to the \cite{halltajv00}
estimates is $62\%$. Regarding the $\alpha=0.5$ case, we found better
performance of the constrained estimators for $V_{123}$, although
lower in magnitude than the $\alpha=0.8$ case, and similar performance
for the bivariate exponent measures. This feature is also supported by
the smaller percentage of change in estimates which is $30\%$. For the
case of strong dependence, i.e., $\alpha=0.2$, the percentage of
change in estimates is very low and equal to $6\%$ which results in
similar efficiency of the estimators for all exponent measures as is
also shown from Figure~\ref{fig:ratio_rmse}.
\begin{table}[htbp!]
  \caption{Monte Carlo estimates of $\tilde{T}_B$ and $\tilde{T}_{B}^c$, $B \in C_3$, 
    for the extreme value logistic case with $\alpha=0.2,0.5$ and $0.8$. }    \vspace{5pt}
  \centering
  \begin{tabular}{ c| c c| c c| c c  }    
    \hline
    &\multicolumn{2}{c|}{$\alpha = 0.2$} & \multicolumn{2}{c|}{ $\alpha = 0.5$}  & \multicolumn{2}{c}{$\alpha = 0.8$}  \\
    $B$& $\tilde{T}_B$  & $\tilde{T}_B^c$& $\tilde{T}_B$  & $\tilde{T}_B^c$ & $\tilde{T}_B$  & $\tilde{T}_B^c$\\
    \hline
    $\{1,2\}$ & 0.02  & 0.02  & 0.36  & 0.35   & 0.76  & 0.72\\
    $\{1,3\}$ & 0.02  & 0.02  & 0.34  & 0.34   & 0.76  & 0.74\\
    $\{2,3\}$ & 0.02  & 0.02  & 0.37  & 0.36   & 0.74  & 0.72\\
    $\{1,2,3\}$& 0.66  & 0.66  & 11.15 & 10.63  & 26.80  & 22.14 
    \label{table:TtildeB}
  \end{tabular}
\end{table}

\noindent Table~\ref{table:TtildeB} shows the Monte Carlo estimates of
the integrated square deviation of the estimators from the theoretical
function. For the case of strong dependence there is no practical
benefit of $\tilde{V}_B^c$ over $\tilde{V}_B$.\ However, in all other
cases the constrained estimators are more efficient than the
unconstrained estimators. This shows that not only does the imposition
of the constraints improve the performance of the estimators for the
higher order exponent measures, but so does for the bivariate level of
dependence. \newline

\noindent To conclude, the performance of the estimators $\tilde{V}_B$
and $\tilde{V}_B^c$ is similar as the dependence increases and becomes
identical in the limiting case of $\alpha \rightarrow 0$. This feature
is explained by the increase in performance of the \cite{halltajv00}
estimators $\tilde{V}_B$ as dependence increases which yields a
consistent set of estimated exponent measures. Overall, we found the
imposition of the new constraints to be beneficial for the simplest
max-stable distribution, i.e., the extreme value logistic, and
superior in efficiency, especially for the case of moderate or weak
dependence. The largest improvement is observed for higher order
exponent measures which is promising for implementations in higher
than 3 dimensions.

\subsection*{Acknowledgements}
I. Papastathopoulos acknowledges financial support from AstraZeneca
and EPSRC.

\end{document}